\documentclass[11pt]{article}
\title{Fibonacci $q$-Gauss sequences}
\author{A.K.Kwa\'sniewski\\  
\\Institute of Computer Science, Bia{\l}ystok University\\
PL-15-887 Bia{\l}ystok, ul.Sosnowa 64, POLAND\\
\\e-mail: kwandr@uwb.edu.pl
\\ Advanced Studies in Contemporary Mathematics vol. 8 (2004) No2 pp.121-124 }

\usepackage{amsmath,amsthm,graphicx}

\chardef\bslash=`\\ 
\newtheorem{hi}{Hint}
\newtheorem{rem}{Remark}

\begin{document}
\maketitle
\begin{abstract}
The summation formula within Pascal triangle resulting in the
Fibonacci sequence is extended to the $q$-binomial coefficients
$q$-Gaussian triangles \cite{1,2}.
\end{abstract}
AMS Classification Numbers: 11B39, 11A39.
\section{Pisa historical remark}
The Fibonacci sequence origin is attributed and referred to the
first edition (lost) of "Fiber abaci" (1202) by Leonardo Fibonacci
[Pisano] (see second edition from 1228 reproduced as Il Liber
Abaci di Leonardo Pisano publicato secondo la lezione Codice
Maglibeciano by Baldassarre Boncompagni in Scritti di Leonardo
Pisano vol. 1, (1857), Rome).
\section{Fibonacci $q$-Gauss sequences}
This note is inspired by \cite{1,2} and by primary binomial
summation formula (\ref{summ1}) for $F_n$:
\begin{equation}\label{summ1}
F_{n+1}=\sum_{k\leq n} \binom{n-k}{k},\quad n\geq 0.
\end{equation}
It is the famous and popular now identity (\ref{summ2}) propagated
via internet \cite{3,4}. Here it is:
\begin{equation}\label{summ2}
F_{n}=\sum_{k=1}^n \binom{n-k}{k-1},\quad n > 0.
\end{equation}
This respectable relation and the property of Pascal triangle is
illustrated by the Figure 1.

\begin{center}
\includegraphics[width=80mm,height=80mm]{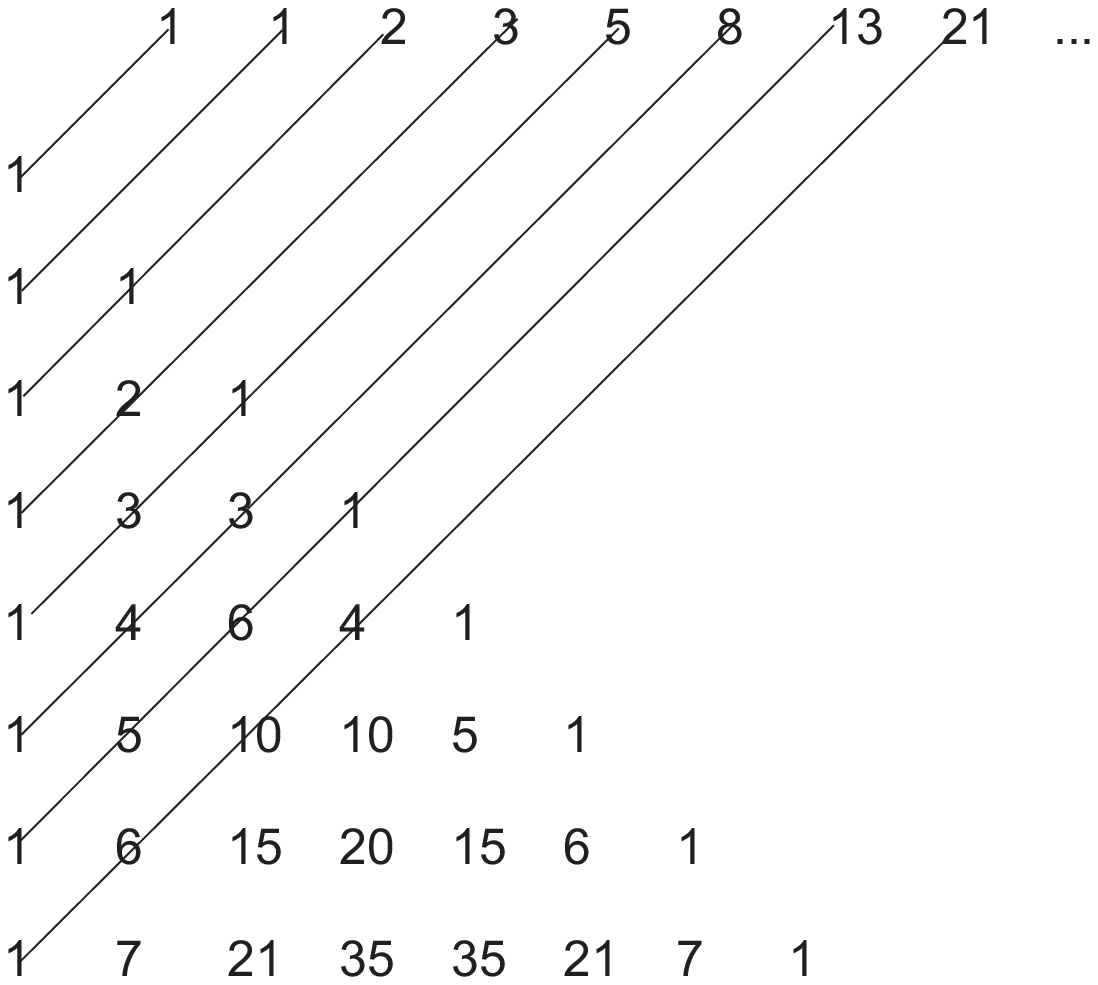}
\end{center}
{\bf Figure 1:} The Fibonacci sequence via summation formula of
Pascal triangle elements.

Formulas (\ref{summ1}) and (\ref{summ2}) are also an immediate
consequence of the $q=1$ recurrence relation:
\begin{equation}\label{contriad}
\left. {{\begin{array}{*{20}l}
{\binom{n+1}{k}_q=\binom{n}{k-1}_q + q^{k}\binom{n}{k}_q} \\
{\binom{0}{0}_q=1}\\
{\binom{k}{0}_q=0\; \textrm{for}\; k>0}
\end{array}} } \right.
\end{equation}
where $n_q!=n_q(n-1)_q!$, $1_q!=0_q!=1$,
$n_q^{\underline{k}}=n_q(n-1)_q(n-2)_q\cdots (n-k+1)_q$,
$\binom{n}{k}_q=\frac{n_q^{\underline{k}}}{k_q!}$.
\begin{rem} {\em (see \cite{1,2})\\
The dual \cite{1} to (\ref{contriad}) recurrence is given by:
\begin{equation}\label{polytriad}
\left. {{\begin{array}{*{20}l}
{x\phi_{k}(x)=q^k \phi_k(x) +\phi_{k+1}(x); \textrm{for}\; k \geq 1}\\
{\phi_0(x)=1}\\
{\phi_{-1}(x)=0}
\end{array}} } \right.
\end{equation}
in accordance with a well known fact that
\begin{equation}\label{extriad}
x^{n}=\sum_{k=0}^{n}\binom{n}{k}_q\phi_k(x),
\end{equation}
$$\phi_{k}(x)=\prod_{i=0}^{k-1}(x-q^i).$$
where $\phi_{k}(x)=\prod_{i=0}^{k-1}(x-q^i)$ are the $q$-Gaussian
polynomials.}
\end{rem}
\begin{rem}{\em (see: \cite{1,2})\\
The connection constants $c_{n,k}=\binom{n}{k}_q$ from relation
(\ref{contriad}) of the triad \cite{1,2} under consideration are
interpreted \cite{5} as the number of $k$-dimensional subspaces in
$n$-th dimensional space $V(n,q)$ over Galois field $GF(q)$.
Restarting this in incidence algebras language \cite{5} we
conclude that in the lattice $L(n,q)$ of all subspaces of
$V(n,q)$:
$$
n_q!=the\;number\;of\;maximal\;chains\;in\;an\;of\;interval\;of\;lenth\;n,
$$
$$
\binom{n}{k}_q= \Big[{\alpha \atop \beta
\gamma}\Big]:=the\;number\;of\;distinct\;elements\;\textrm{z}\;in\;a\;segment\;
[\textrm{x,y}]
$$
$$of\;type\;\alpha\;such\;that\;[\textrm{x,z}]\;is\;of\;type\;\beta\;while\;
[\textrm{z,y}]\;is\;of\;type\;\gamma,
$$
where $\textrm{x,y,z} \in L(n,q)$ and the type of $[\textrm{x,y}]$
is equal to dim[x/y]. (Here [y/x] denotes the subtraction of
linear subspaces \cite{5}).}
\end{rem}
The corresponding to (\ref{contriad}) $q$-Gaussian-Pascal
triangles have analogous summation formula property as the one
given by (\ref{summ1}) and illustrated by the Figure 1. Let it be
now illustrated for the case $q=2$ by the Figure 2 showing the
first elements of the sequence $\{F_n^{[q^j]}\}$ ($j=0$, $q=2$):
\begin{equation}\label{F-G}
\left. {{\begin{array}{*{20}l}
{F_{n+1}^{[q^j]}=\sum \limits_{k\leq n} \binom{n-k}{k}_q(q^j)^{k-1},\quad n,j\geq 0,\;\;q\geq 1}\\
{F_0^{[q^j]}=0}
\end{array}} } \right.
\end{equation}

\begin{center}
\includegraphics[width=80mm,height=80mm]{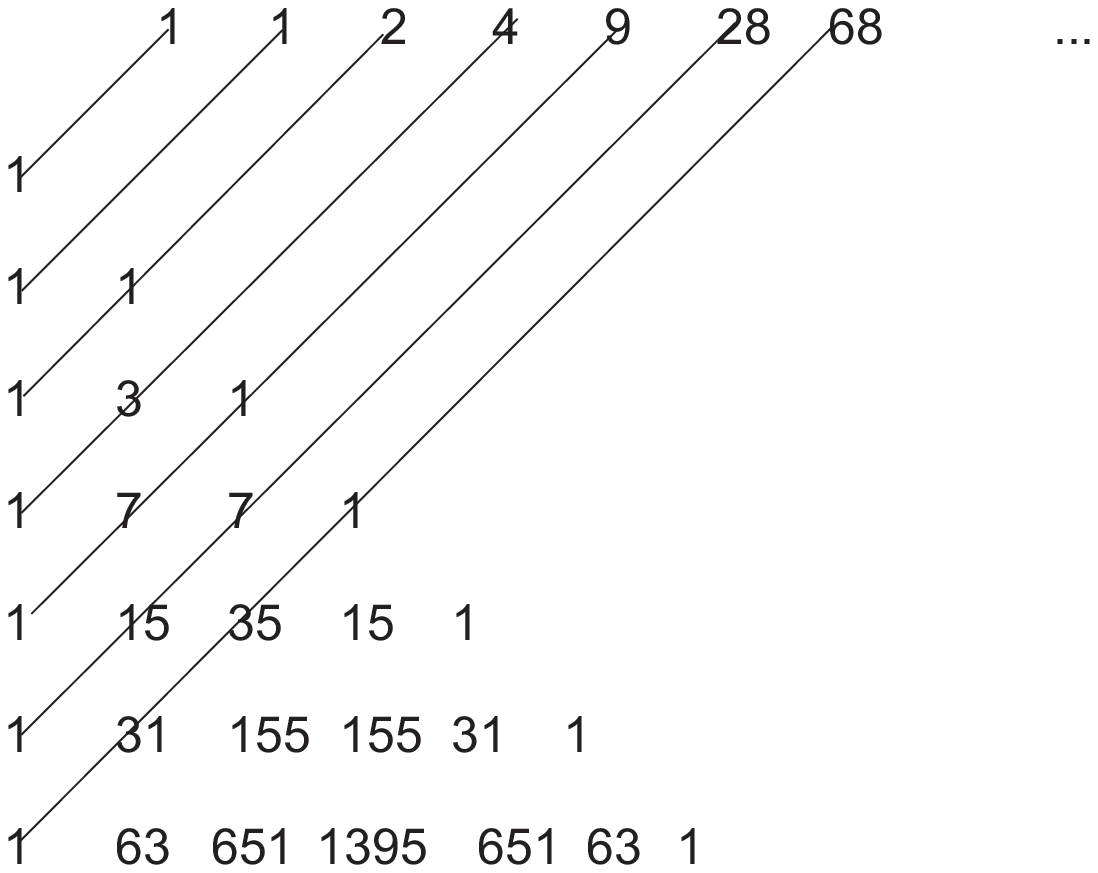}
\end{center}
{\bf Figure 2:} The Fibonacci $q$-Gauss sequence via summation
formula of $q$-Gauss Pascal triangle elements; $q=2$, $j=0$.

We shall call the sequence of sequences defined by (\ref{F-G}) and
with given "the Fibonacci $q$-Gauss family" of $j$-Pascal
$q$-Gauss sequences. For $q=1$ (\ref{F-G}) becomes (\ref{summ2})
and we get the Pascal triangle case with the resulting Fibonacci
sequence via (\ref{summ1}) or (\ref{summ2}).

Using the recurrence (\ref{contriad}) one derives - in a way one
does that in $q=1$ case - the following recurrence relations for
the Fibonacci $q$-Gauss family sequences:
\begin{equation}\label{recF-G}
\left. {{\begin{array}{*{20}l}
{F_{n+2}^{[q^j]}=F_{n+1}^{[q^{j+1}]}+F_{n}^{[q^j]},\quad n,j \geq
0,\;q\geq 1,}\\
{F_0^{[q^j]}=0}\\
{F_1^{[q^j]}=q_0^j}
\end{array}} } \right.
\end{equation}
\begin{hi}{\em One may extend the notion of the Fibonacci
$q$-Gauss family to the notion of the Fibonacci triad triangle
\cite{1,2} family with the corresponding to (\ref{recF-G}) "stream
of recurrences".}
\end{hi}
\begin{hi}{\em The ordinary generating function for the Fibonacci
$q$-Gauss family given by
\begin{equation}
\left. {{\begin{array}{*{20}l} {F(q;z,x)=\sum \limits_{l,m \geq
0}F(q^l;x)z^m}\\
{F(q^l;x)=\sum \limits_{n\geq 0}F_n^{[q^j]}x^n}
\end{array}} } \right.
\end{equation}
becomes the well known one in the case of $q=1$:
\begin{equation}
F(1;0,x)=\frac{x}{1-x-x^2}.
\end{equation}}
\end{hi}
Announcement on an investigation of the properties of the ordinary
generating function for the Fibonacci $q$-Gauss family we postpone
for the subsequent report on the subject.

\end{document}